\documentclass{amsart}
\usepackage{amsmath,amssymb,amsxtra,amsthm,amscd}
\usepackage[mathscr]{eucal}
\usepackage{bbm}
\usepackage[all,cmtip]{xy}
\usepackage{lmodern}
\usepackage{graphicx}
\usepackage{todonotes}
\usepackage{algorithmic}
\usepackage{algorithm}
\usepackage{eqparbox}
\usepackage{clrscode}
\usepackage{multirow}
\usepackage{epic}

\usepackage{tikz} 

\newif\ifShowLabels
\ShowLabelstrue
\newcommand{\TeXref}[1]{
\marginpar{\scriptsize \texttt{#1}}}
\ShowLabelsfalse

\DeclareMathOperator{\J}{\mathit{J}}

\DeclareMathOperator*{\one}{1}
\newcommand{\onehatplace}[1]
{ \one^{\substack{#1 \\ \frown}} }

\DeclareMathOperator*{\bones}{\times}
\newcommand{\undertimes}[1]
{ \bones_{#1} }

\DeclareMathOperator*{\bowl}{\cup}
\newcommand{\undercup}[1]
{ \bowl_{#1} }

\DeclareMathOperator*{\arch}{\cap}
\newcommand{\undercap}[1]
{ \arch_{#1} }

\newcommand{\pull}
{\!\!\! -\!\!\! -\!\!\! -\!\!\!}

\DeclareMathOperator*{\holimprep}{holim}                       
\newcommand{\holim}[1]%
{\displaystyle\holimprep_{\substack{\leftarrow \pull - \\ #1}} \, }

\DeclareMathOperator*{\hocolimprep}{hocolim}                   
\newcommand{\hocolim}[1]%
{\displaystyle\hocolimprep_{\substack{- \pull \rightarrow \\ #1}} \, }

\DeclareMathOperator*{\plainlim}{lim}                           
\newcommand{\contralim}[1]%
{\displaystyle\plainlim_{\substack{\leftarrow \pull - \\ #1}} \, }

\DeclareMathOperator*{\plaincolim}{colim}                       
\newcommand{\colim}[1]%
{\displaystyle\plaincolim_{\substack{- \pull \rightarrow \\ #1}} \, }

\DeclareMathOperator*{\laxlimplain}{laxlim}                     
\newcommand{\laxlim}[1]%
{\displaystyle\laxlimplain_{\substack{\leftarrow \pull - \\ #1}} \, }

\newcommand\LONGCOMMENT[1]{%
  \hfill\#\ \begin{minipage}[t]{\eqboxwidth{COMMENT}}#1\strut\end{minipage}%
}

\theoremstyle{plain}
\newtheorem{Thm}{Theorem}[section]

\newtheorem{Cor}[Thm]{Corollary}

\newtheorem{Lem}[Thm]{Lemma}
\newtheorem{Prop}[Thm]{Proposition}

\theoremstyle{definition}
\newtheorem{Def}[Thm]{Definition}

\newtheorem{Ex}[Thm]{Example}

\newtheorem{Rem}[Thm]{Remark}

\theoremstyle{remark}
\newtheorem{Not}[Thm]{Notation}

\newtheoremstyle{freestylethm}{6pt}{6pt}{\itshape}{}%
                {\bfseries}{}{.5em}{\thmnote{#3}}
\theoremstyle{freestylethm}

\newcommand{\SecRef}[2]{\section{#1}\label{S:#2}%
\ifShowLabels \TeXref{{S:#2}} \fi}

\newcommand{\refS}[1]{\textup{\ref{S:#1}}}

\newcommand{\refT}[1]{\textup{\ref{T:#1}}}

\newcommand{\refD}[1]{\textup{\ref{D:#1}}}
\newcommand{\refC}[1]{\textup{\ref{C:#1}}}
\newcommand{\refE}[1]{\textup{\ref{E:#1}}}

\newcommand{\refR}[1]{\textup{\ref{R:#1}}}

\newenvironment{ThmRef}[1]%
{ \begin{Thm} \label{T:#1}
\ifShowLabels \TeXref{T:#1} \fi }%
{ \end{Thm} }
\newenvironment{DefRef}[1]%
{ \begin{Def} \label{D:#1}
\ifShowLabels \TeXref{D:#1} \fi }%
{ \end{Def} }
{ \begin{Lem} \label{L:#1}
\ifShowLabels \TeXref{L:#1} \fi }%
{ \end{Lem} }
\newenvironment{CorRef}[1]%
{ \begin{Cor} \label{C:#1}
\ifShowLabels \TeXref{C:#1} \fi }%
{ \end{Cor} }
\newenvironment{RemRef}[1]%
{ \begin{Rem} \label{R:#1}
\ifShowLabels \TeXref{R:#1} \fi }%
{ \end{Rem} }
\newenvironment{PropRef}[1]%
{ \begin{Prop} \label{P:#1}
\ifShowLabels \TeXref{P:#1} \fi }%
{ \end{Prop} }
{ \begin{Ex} \label{E:#1}
\ifShowLabels \TeXref{E:#1} \fi  }%
{ \end{Ex} }
\newenvironment{NotRef}[1]%
{ \begin{Not} \label{N:#1}
\ifShowLabels \TeXref{N:#1} \fi }%
{ \end{Not} }

{ \begin{Thm} [#2]\label{T:#1}
\ifShowLabels \TeXref{T:#1} \fi }%
{ \end{Thm} }
{ \begin{Def} [#2]\label{D:#1}
\ifShowLabels \TeXref{D:#1} \fi }%
{ \end{Def} }
{ \begin{Lem} [#2]\label{L:#1}
\ifShowLabels \TeXref{L:#1} \fi }%
{ \end{Lem} }
{ \begin{Cor} [#2]\label{C:#1}
\ifShowLabels \TeXref{C:#1} \fi }
{ \end{Cor} }
{ \begin{Rem} [#2]\label{R:#1}
\ifShowLabels \TeXref{R:#1} \fi }%
{ \end{Rem} }
{ \begin{Prop} [#2]\label{P:#1}
\ifShowLabels \TeXref{P:#1} \fi }%
{ \end{Prop} }
\newenvironment{ExRefName}[2]%
{ \begin{Ex} [#2]\label{E:#1}
\ifShowLabels \TeXref{E:#1} \fi }%
{ \end{Ex} }

\begin{document}

\title[Merging discrete Morse Vector fields]{Merging discrete Morse vector fields:\\ a case of stubborn geometric parallelization}
\author{Douglas Lenseth and Boris Goldfarb}
\address{Department of Mathematics and Statistics, State University of New York, Albany, NY 12222}
\email{dlenseth@albany.edu, bgoldfarb@albany.edu}



\date{\today}

\begin{abstract}
 We address the basic question in discrete Morse theory of combining discrete gradient fields that are partially defined on subsets of the given complex.
This is a well-posed question when the discrete gradient field $V$ is generated using a fixed algorithm which has a local nature.  One example is \texttt{ProcessLowerStars}, a widely used algorithm for computing persistent homology associated to a grey-scale image in 2D or 3D.  While the algorithm for $V$ may be inherently local, being computed within stars of vertices and so embarrassingly parallelizable, in practical use it is natural to want to distribute the computation over patches $P_{i}$, apply the chosen algorithm to compute the fields $V_{i}$ associated to each patch, and then assemble the ambient field $V$ from these.  Simply merging the fields from the patches, even when that makes sense, gives a wrong answer.  
We develop both very general merging procedures and leaner versions designed for specific, easy to arrange covering patterns.
\end{abstract}

\maketitle

\tableofcontents

\SecRef{Introduction}{MT}
\par Discrete Morse theory \cite{F:02} is a fairly new and powerful tool, created by Robin Forman in the 1990s, that has many applications in many fields.  One important application of discrete Morse theory is to streamlining homology calculations.  It has been a very useful tool in topological data analysis, for example in computing persistent homology. 

\par More specifically, there are algorithms that take both 2D and 3D digital grayscale images and create discrete Morse functions from their grayscale function values defined on pixels or voxels. These images give rise to cubical complexes, where the voxels or pixels are the 0-cells of the complex. One algorithm, \texttt{ProcessLowerStars} from \cite{RWS:11}, takes all of the cells in the cubical complex and puts them either singly in the set of critical cells, $C$, or as a face-coface (of codimension 1) pair in the discrete vector field, $V$.  The set of critical cells and the discrete vector field are the defining features of a discrete Morse function, and in turn the discrete Morse complex.  In discrete Morse theory, the discrete vector field determines the maps in the simplified chain complex of the Morse complex that is created from the critical cells of the discrete Morse function.

Given a general algorithm $\alpha$ applied to cells in a regular CW complex $K$ which takes data from $k$-neighborhoods of cells in $K$ and uses that data to either pair cells together in a list which is to become the discrete vector field $V$ or place the cells singly in the list of critical cells $C$.
The basic question we address is how do we merge the vector fields in the elementary situation where $K=U \cup W$, and $U$, $W$ are two subcomplexes to which $\alpha$ is applied individually to get the discrete vector fields $V(U)$ and $V(W)$? The goal is to reconstruct the correct vector field $V(K)$ that we would get on all of $K$.  It turns out this is not straightforward, in that $V(K) \neq V(U) \cup V(W)$, and so $\alpha$ is not embarassingly parallelizable.  Being able to do such a computation will involve looking at where mistakes occur in naively merging the vector fields together as in the right-hand side of the formula (or, in general, related vector fields) and ``crossing out'' those mistakes.

\par  What do we mean by an algorithm being embarassingly parallelizable?  One can take an example with operators in calculus to illustrate the difference between embarrassingly and stubbornly parallelizable algorithms. Differentiation is embarassingly parallelizable with respect to addition and subtraction, i.e. $(f(x) \pm g(x))'=f'(x) \pm g'(x)$.  But it is stubbornly parallelizable with respect to multiplication and division, i.e. $(f(x)g(x))'=f'(x)g(x)+f(x)g'(x)$, which is a non-obvious replacement of the straightforward but incorrect formula $(f(x)g(x))'=f'(x)g'(x)$.  

We will present the main parallelization theorem, Theorem \refT{NEJCHCB}, that will allow us to do this merging for a general uniformly $k$-local algorithm $\alpha$.  Here $k$-local refers to the nature of the algorithm in that it makes the classification decision about a cell or a pair by processing only the information from $k$-neighborhoods of cells in a regular CW complex.  This turns out to be the common type of algorithms in practice.  The authors don't know any $\alpha$ that builds discrete Morse fields or functions and doesn't have this property.  We then make this merging process more efficient in Corollary \refC{MNIH} by noting some geometric properties of the vector fields we get, which leads to a merging that involves applying $\alpha$ on fewer subcomplexes, i.e. a more efficient distributed computation.  We then look at the situation where the subcomplexes $U$ and $W$ are in antithetic position, that is satisfying $(U \cap W)[k]=U[k]\cap W[k]$, where $k \geq 1$ and $U[k]$ is the $k$-neighborhood of the subcomplex $U$.  This special relative positioning of the patches allows to further improve the efficiency by allowing to use smaller auxiliary patches.  This is spelled out in Corollary \refC{NEJCHCB}.

\par In the context of $\alpha$ being \texttt{ProcessLowerStars}, it turns out the rigid structure of the cubical complexes that are formed from the images gives a very efficient formula, related to but not exactly matching the formula for the more general case.  We will present a formula, in the 2D image case, that allows us to take lists of discrete vector fields directly from two patches and merge them together to get the discrete vector field of the union of the two patches, under the condition that the two patches must have an overlap that contains at least one $2$-cell.  This is done in section \refS{PA} and is an illustration of how the general formulas can be further tweaked to improve the efficiency of the procedure.

\par Why do we want to parallelize this \texttt{ProcessLowerStars} algorithm?  For one, if one has a particularly large image with a lot of pixels or voxels, being able to break the image up and apply the algorithm on each piece will reduce the time cost compared to applying the algorithm on the whole image.  It also may be that an image may come to you in pieces, \'a la online machine learning.  The methods of this paper allow to process each piece as it arrives and either gradually build the discrete field or put the whole field together as a final quick step.

\par  In an application to the explicit assembly strategy for the \texttt{ProcessLowerStars} algorithm on a 2D image, we will assume that our computer has a hypothetical limited constraint (i.e. restricted to processing at most a certain number of pixels).  In this situation we will have to break up the image so that the number of pixels in each patch is less than the constraint of our computer.  These patches will be rectangular and will allow us to index our patches both in the horizontal and vertical directions.  We will apply the \texttt{ProcessLowerStars} algorithm to each patch (and all intersections, in both the horizontal and vertical direction, of adjacent patches).  We create an algorithm that merges the discrete vector field in each of the patches to obtain the discrete vector field of the whole image, using the formula on two patches.  In order to do this, our algorithm will merge starting in the top left corner of the image and move along each horizontal strip and their intersections.  Once the discrete vector field of every horizontal strip and their intersections is found, the algorithm then moves vertically from top down to merge the discrete vector fields of the horizontal strips two at a time again making use of the formula that will be provided.

\par Finally we will look at another specific special geometric situation that involves simplicial trees in section \refS{EX}.  This is done with two purposes.  It is known that finite products of trees contain coarse images of the most common structures in geometry of finite asymptotic dimension.  A sufficient partition scheme in a product of trees will therefore provide a generic template to approach most common regular cellular complexes.  On the other hand, we will find that a tree is an example that distinguishes the advantages of each of the two main parallelization results of the paper.

We finish with a collection of remarks on applications of the results, the literature, and directions for future work in sections \refS{D} and  \refS{Disc}.

\SecRef{Discrete gradient vector fields}{BGOAG}

All complexes in this paper are regular CW complexes, in the sense that characteristic maps of all cells are homeomorphisms onto the closures of cells in the complex.  We refer for the general background in discrete Morse theory on regular CW complexes to \cite{Knudson}.

\begin{NotRef}{QWRC}
	Suppose $K$ is a regular complex, and suppose $f$ is a discrete Morse function on $K$. Then there are well-defined sets of critical cells $C (K,f)$ and the discrete vector field that can be thought of as a set of pairs of regular cells $V (K,f)$ or $V (K)$. When we use the notation $(\sigma \le \tau )$ for a pair in $V (K)$, we mean $\sigma$ is a codimension 1 face of $\tau$.
\end{NotRef}

\begin{DefRef}{Adj}
We introduce the following relation among cells in $K$.  Two cells $\sigma_1$ and $\sigma_2$ are \textit{adjacent} if the closures of the two cells in $K$ have a point in common.	The \textit{star} of a cell $\sigma$ is the smallest cellular subcomplex of $K$ that contains all cells adjacent to $\sigma$.  Given a collection of cells $\mathcal{C}$ of $K$, let the \textit{star} of $\mathcal{C}$ be the smallest subcomplex that contains the stars of all of the cells in $\mathcal{C}$. If $K'$ is a subcomplex, we also call its star the \textit{1-neighborhood} of $K'$ denoted $K'[1]$.  The star of a star is the 2-neighborhood, etc.  We will use the notation $K' [n]$ for the $n$-neighborhood of $K'$.  
The $k$-\textit{border} of a subcomplex $L$ of $K$ consists of those cells in $L$ whose $k$-neighborhoods are not contained entirely in $L$. 
\end{DefRef}

We will assume that the lists $C (K,f)$ and $V (K,f)$ are results of applying an algorithm $\alpha \colon K \mapsto V$  of the following nature.  

\begin{DefRef}{klocal}
Given a face-coface pair of cells $(\sigma \le \tau)$, the algorithm decides whether the pair $(\sigma \le \tau)$ is placed in the list $V (K)$ based on data specific to $\alpha$ coming from some $k$-neighborhood $\sigma[k]$, for some uniformly fixed number $k > 0$.  The result of the algorithm is a discrete vector field $V_{\alpha} (K)$ which can be therefore realized as $V (K,f)$ for some discrete Morse function $f$.  We will say $\alpha$ with this property is \textit{uniformly local} or, more precisely, \textit{uniformly $k$-local}.
\end{DefRef}

In fact, in all applications that generate $V(K)$ the authors see in the literature, it is always done using a uniformly $1$-local algorithm.  This isn't surprising of course, as traditional smooth vector fields are generated through differentiation procedures which are inherently local.  

\begin{ExRefName}{mki}{\texttt{ProcessLowerStars} is uniformly local}
 An example of such algorithm is \texttt{ProcessLowerStars} which produces $V(K)$ starting with an arbitrary positive bounded function (a.k.a grayscale function) on cells (pixels or voxels) of a 2D or 3D pixel grid $K$.  Since the algorithm qualifies $(\sigma \le \tau)$ based on the values of a function on vertices in $\sigma [1]$, this algorithm is uniformly $k$-local for $k=1$. 
 
 To be more specific, to make these images into a regular CW complex, we will translate them to a cubical complexes by making the pixels or voxels correspond to 0-cells of the cubical complex, $K$.  This gives us a positive bounded function, call it $g$, on a cubical complex, more specifically, a cubulated plane or space.  To apply the algorithm, it is ideal to have the values of $g$ on all of the 0-cells be unique.  If they are not, one can make small changes (i.e., a linear ramp) to $g$ to ensure they are unique.  Once there are unique values of $g$, the algorithm inspects all of the cells in the lower star of a 0-cell, $x$. The lower star of $x$, $L(x)$, contains all cells $\sigma \in K $ in the cubical complex such that $g(x)=\max\limits_{y \in \sigma} g(y)$. To give an ordering to higher dimensional cells, a new function, $G$, is introduced, defined as follows:  If $\sigma$ contains the vertices $\{x,y_{1}, \cdots,y_{n}\}$, then $G(\sigma)=\{g(x),g(y_{i_{1}}), \cdots, g(y_{i_{n}})\}$, where $g(x)>g(y_{i_{1}})>\cdots>g(y_{i_{n}})$. This will allow us to impose the lexicographical ordering on $G$ when performing this algorithm.
 
 The algorithm, itself, will then take all of the cells of in $L(x)$ and either possibly pair a cell, $\tau$ with a codimensional 1 face, $\sigma$ and place them $(\sigma \le \tau)$ in $V$ or take a cell, $\sigma$ and place it singly in $C$ in the following way: 
 
\begin{itemize}
    \item If $L(x)=\{x\}$, then $x \in C$.  Otherwise, take the minimal (with respect to $G$) $1$-cell, $\tau$, and $(x \le \tau) \in V$. All other $1$-cells are added to a queue called PQzero, since they have no remaining unpaired faces in $L(x)$.  All cofaces of codmension 1 of $\tau$ in $L(x)$ that have 1 unpaired face in $L(x)$ are added to a different queue called PQone.
    \item The algorithm then takes the minimal cell (with respect to $G$) that is in PQone and  either moves it to PQzero, if it has no unpaired faces remaining in $L(x)$, or if it still has an unpaired face it gets paired with that face and is put into $V$. Then you look at all cofaces of codimension 1 of both cells that were just paired and put into $V$.  If any of these cofaces have exactly one unpaired face, they are added to the queue PQone.
    \item  Then, if PQzero is not empty, it takes the minimal cell (with respect to $G$), call it $\delta$, and places it singly in $C$. Then all cofaces of codimension 1 of $\delta$ are inspected.  If any these cofaces have exactly one unpaired face, it is placed in PQone. 
\end{itemize}
 This will keep going until both PQzero and PQone are empty.  Note that, in particular $L(x)=x[1]$, using the notation defined above. The \texttt{ProcessLowerStars}, in fact, only looks at $\sigma[1]$, where $\sigma$ is a $0$-cell.
\end{ExRefName}

\begin{ExRefName}{YIH}{failure of naive merge}
	We present a simple example with $V_{\alpha} (K) \ne V_{\alpha} (U) \cup V_{\alpha} (W)$, where $\alpha$ is the \texttt{ProcessLowerStars} algorithm from \cite{RWS:11} and $U$ and $W$ are subcomplexes of $K$ such that $K=U\cup W$. 
Let $K$ be the following 1-dimensional cubical complex with the function values, $f$, on the $0$-cells.

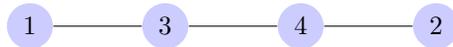
\begin{figure}[h!]
	\centering
\begin{tikzpicture}  
  [scale=.9,auto=center,every node/.style={circle,fill=blue!20}] 
  \node (a1) at (1,0) {1};  
  \node (a2) at (3,0) {3}; 
  \node (a3) at (5,0) {4};  
  \node (a4) at (7,0) {2};  
  \draw (a1) -- (a2); 
  \draw (a2) -- (a3);  
  \draw (a3) -- (a4);
\end{tikzpicture}  
	\caption{1-D example of a failed merging of discrete vector fields.}
	\label{Merge}
\end{figure}

First, we apply $\alpha$ to all of $K$. This gives $V_{\alpha}(K)=\{\{3, 31\}, \{4, 42\}\}$.  Next, we will break up $K$ into the following subcomplexes: $U=\{1,3,4,31,43\}$ and $W=\{2,3,4,43,42\}$ .  When we apply $\alpha$ to $U$, we get $V_{\alpha}(U)=\{\{3,31\},\{4,43\}\}$.  When we apply $\alpha$ to $W$, we get $V_{\alpha}(W)=\{\{4,42\}\}$. Clearly, $V_{\alpha} (K) \ne V_{\alpha} (U) \cup V_{\alpha} (W)$, so simply merging the lists from the subcomplexes does not give us the correct discrete vector field for all of $K$.  In this particular instance, we get an extra pairing, the $\{4,43\}$ from $U$. 
\end{ExRefName}

The failure of naive merging is made clearer using the following definition and explanation.  Let $\sigma_C [k]$ be the $k$-star of $\sigma$ viewed as a cell within the subcomplex $C$ of $K$.  This is an intrinsic to $C$ construction and so may differ from $\sigma [k]$, which is the same as $\sigma_K [k]$, for some combinations of $\sigma$ and $k$.  Since the algorithm $\alpha$ takes as input the data from $\sigma_C [k]$, there are cases when the discrepancy leads to different constructions of the vector fields as Example \refE{YIH} shows.  This also points to some cases when this discrepancy does not happen.

\begin{PropRef}{BQA}
Suppose $K$ is a disjoint union of subcomplexes $C_i$, for $i$ from some index set $I$.  In other words, a connected subset of $K$ is contained entirely within one and only one $C_i$. Then $V_{\alpha} (K) = \bigcup_{i \in I} V_{\alpha} (C_i)$.
\end{PropRef}

\begin{proof}
    For every pair $(\sigma \le \tau)$ its inclusion in $V_{\alpha} (C_i)$ is decided by $\alpha$ from the data within $\sigma_{C_i} [k]$.  Since in our case $\sigma_{C_i} [k] = \sigma_K [k]$, that decision is precisely the same within $C_i$ as within $K$.
\end{proof}

This argument makes it clear that failure of the simple merge is caused by possible discrepancies between $\sigma_{C} [k]$ and $\sigma_K [k]$ for cells $\sigma$ in the $k$-border of $C$ (see Definition \refD{Adj}).

\SecRef{The parallelization theorems}{FT}

Suppose a complex $K$ is assigned a gradient vector field $V$ using an algorithm $\alpha$. 
Suppose further that $U$ and $W$ are two subcomplexes of $K$ that form a covering of $K$.  In particular, there are vector fields $V_{\alpha} (K)$, $V_{\alpha} (U[i])$, and $V_{\alpha} (W[i])$ for various values of $i$.

\begin{ThmRef}{NEJCHCB}
If  $\alpha$ is a uniformly $k$-local algorithm for some $k \ge 1$, then
  \begin{multline*}
  	V_{\alpha} (K) = \left( V_{\alpha} (U[k]) \cup V_{\alpha} (W[k]) \right) \setminus \notag \\
  	 (V_{\alpha} (U[k]) \cup V_{\alpha} (W[k]) \setminus V_{\alpha} (U[2k+1] \cap W[2k+1]). \notag
  \end{multline*}
In slightly different terms,
  \begin{multline*}
  	V_{\alpha} (K) = \left( V_{\alpha} (U[k]) \cup V_{\alpha} (W[k]) \right) \setminus \notag \\
  	\left( V_{\alpha} (U[k]) \setminus V_{\alpha} (U[2k+1]) \right) \setminus
  \left( V_{\alpha} (W[k]) \setminus V_{\alpha} (W[2k+1]) \right). \notag
  \end{multline*}
\end{ThmRef}

The informal idea is this.  It will be elementary to check that $V_{\alpha} (K)$ is contained in $V_{\alpha} (U[k]) \cup V_{\alpha} (W[k])$ using that the algorithm $\alpha$ is uniformly $k$-local. Let's refer to the difference as ``mistakes'' in building $V_{\alpha} (K)$ from the union.  We will see that for the same reason, $V_{\alpha} (U[k]) \setminus V_{\alpha} (U[2k+1])$ can be viewed as all the mistakes made near the border of $U[k]$.  This set is interpreted as pairs of cells in $V_{\alpha} (U[k])$ that are not in $V_{\alpha} (K)$.  Similarly, $V_{\alpha} (W[k]) \setminus V_{\alpha} (W[2k+1])$ is the set of mistakes made near the border of $W[k]$. Since the cells near the borders of the subcomplexes are far enough away in the adjacency relation chains, they do not contain common adjacent cells. The local nature of the algorithm, again, guarantees that no other mistakes are made away from the borders of $U[k]$ and $W[k]$.

We now give a formal proof.
	
\begin{proof}
	If the pair $(\sigma \le \tau)$ is contained in $V_{\alpha}(K)$, it was placed in this list according to the algorithm $\alpha$.  According to the assumption, this placement is the outcome of inspection of data relevant to $\alpha$ in $\sigma [k]$.  If $\sigma$ and $\tau$ are cells of $U$ then $\sigma [k]$ is contained in $U[k]$.  Since $\alpha$ applied to $U[k]$ would use the exact same data as $\sigma [k] = \sigma_{U[k]} [k]$, this proves that $(\sigma \le \tau)$ is in $V_{\alpha}(U[k])$.  Repeating the argument for $W$ in place of $U$,  this gives the inclusion $V_{\alpha} (K) \subset V_{\alpha} (U[k]) \cup V_{\alpha} (W[k])$.
	
	We now examine $(\sigma \le \tau)$ in the difference
	\[
	M = V_{\alpha} (U[k]) \cup V_{\alpha} (W[k]) \setminus V_{\alpha} (K).
	\]
	
	Observe that whenever the cell $\sigma$ is outside of $W[k]$, it is contained in $U$.  This makes $\sigma_K [k] = \sigma_{U[k]} [k]$, so the outcome of the decision $\alpha$ makes on inclusion of $(\sigma \le \tau)$ in respectively $V_{\alpha} (K)$ and $V_{\alpha} (U[k])$ is the same.  This shows that such $(\sigma \le \tau)$ is disqualified from $M$ as it is excluded as part of $V_{\alpha} (K)$.  Our conclusion is that whenever $(\sigma \le \tau)$ is in $M$, $\sigma$ must be in $W[k]$. A symmetric argument proves that it also must also be in $U[k]$.  
	
	The last paragraph proved that our pair of cells $\sigma$, $\tau$ is contained in $U[k+1] \cap W[k+1]$.  Therefore $\alpha$ uses the same data for placing $(\sigma \le \tau)$ in either $V_{\alpha} (U[2k+1] \cap W[2k+1])$ or $V_{\alpha} (K)$.
	This allows us to rewrite
	\[
	M = V_{\alpha} (U[k]) \cup V_{\alpha} (W[k]) \setminus V_{\alpha} (U[2k+1] \cap W[2k+1]).
	\]
	This gives the first formula.
	
	Let $M_U = V_{\alpha} (U[k]) \setminus V_{\alpha} (K)$ and $M_W = V_{\alpha} (W[k]) \setminus V_{\alpha} (K)$, so that $M = M_U \cup M_W$. The same argument as above shows that
	\[
	M_U = V_{\alpha} (U[k]) \setminus V_{\alpha} (U[2k+1]),
	\]
	and similarly
	\[
	M_W = V_{\alpha} (W[k]) \setminus V_{\alpha} (W[2k+1]).
	\]
Now 
 \[
  	V_{\alpha} (K) =  V_{\alpha} (U[k]) \cup V_{\alpha} (W[k]) \setminus 
  	\left( M_U \cup M_W \right)
  \]
  gives the second formula.
\end{proof}

\begin{RemRef}{HJQAS}
Our interest in this theorem is that it allows to express the global list $V_{\alpha} (K)$ in terms of generally smaller partial lists $V_{\alpha} (U[k])$, $V_{\alpha} (W[k])$, $V_{\alpha} (U[2k+1])$, and $V_{\alpha} (W[2k+1])$.
\end{RemRef}

The computation of these smaller lists can be further simplified  by the observation that the difference between $V_{\alpha} (U[k])$ and $V_{\alpha} (U[2k+1])$ is contained in $V_{\alpha} (U[2k+1] \cap W[2k+1])$.

\begin{CorRef}{MNIH}
If  $\alpha$ is a uniformly $k$-local algorithm for some $k \ge 1$, then
  \begin{multline*}
  	V_{\alpha} (K) = \left( V_{\alpha} (U[k]) \cup V_{\alpha} (W[k]) \right) \setminus \notag \\
  	\left( V_{\alpha} (U[k] \cap W[k]) \right)
  \setminus
  \left( V_{\alpha} (U[2k+1] \cap W[2k+1]) \right). \notag
  \end{multline*}
\end{CorRef}

\begin{proof}
The formula from Theorem \refT{NEJCHCB} is rewritten in terms of the intersection $U[2k+1] \cap W[2k+1]$.  Notice that the restriction to $V_{\alpha} (U[k] \cap W[k])$ is justified because the $k$-borders of $U$ and $W$ and thus also $M$ are contained in $U[k] \cap W[k]$.
\end{proof}

We next implement the formula in this simple case of a two-set covering.

\begin{algorithm}
\caption{Two-set parallelization from Corollary \refC{MNIH}}
\begin{algorithmic}
\STATE{\bfseries Inputs:} \begin{itemize}
    \item $k \ge 1$
    \item $\alpha$ uniformly $k$-local algorithm
    \item $K$ regular CW complex
    \item $D$ data attached to cells in $K$, specific for use in $\alpha$
    \item $U(1)$, $U(2)$ a covering of $K$ by subcomplexes given as subsets of $K$
    \end{itemize}
\STATE{\bfseries Outputs:} \begin{itemize}
    \item $V = V_{\alpha} (K,D)$ discrete vector field
    \item $C = C_{\alpha} (K,D)$ critical cells
    \end{itemize}
\FOR{$i$ = $1$, $2$} 
    \STATE{} \begin{itemize}
    \item \textbf{build} $U(i)[k]$
    \item \textbf{build} $U(i)[2k+1]$
    \item \textbf{build} $D\vert U(i)[k]$
    \item \textbf{evaluate} $\alpha$: $U(i[k])$, $D\vert U(i)[k] \to A(i,k)$
    \end{itemize} 
    \LONGCOMMENT{preprocessing step to obtain discrete vector fields in enlargements of patches}
\ENDFOR
\STATE{\textbf{build} the union $V(1,2) = A(1,k) \cup A(2,k)$}
\STATE{\textbf{build} the intersection $U(1,2,k) = U(1)[k] \cap U(2)[k]$}
\STATE{\textbf{build} the intersection $U(1,2,2k+1) = U(1)[2k+1] \cap U(2)[2k+1]$}
\STATE{\textbf{build} $D\vert U(1,2,k)$}
\STATE{\textbf{build} $D\vert U(1,2,2k+1)$}
\STATE{\textbf{evaluate} $\alpha$: $U(1,2,k)$, $D\vert U(1,2,k) \to A(1,2,k)$}
\STATE{\textbf{evaluate} $\alpha$: $U(1,2,2k+1)$, $D\vert U(1,2,2k+1) \to A(1,2,2k+1)$}
\STATE{\textbf{build} the difference $D = A(1,2,k) \setminus A(1,2,2k+1)$}
\STATE{\textbf{save} the difference $V = V(1,2) \setminus D$}
\STATE{\textbf{build} $N$ as all cells employed in $V$}
\STATE{\textbf{save} the difference $C = K \setminus N$}
\end{algorithmic}
\end{algorithm}

There is a special geometric situation with a more direct identification of intersections of enlargements. 

\begin{DefRef}{NAT}
Two subcomplexes $U$ and $W$ of $K$ are said to be in \textit{antithetic position} if $(U \cap W) [n] = U[n] \cap W[n]$ for all $n \ge 0$. 
\end{DefRef}

We can now state the following consequence of Corollary \refC{MNIH}.

\begin{CorRef}{NEJCHCB}
Suppose $\alpha$ is an algorithm that is uniformly $k$-local for some $k \ge 1$ and suppose $U$, $W$ is an antithetic pair, then
  \[
  	V_{\alpha} (K) = \left( V_{\alpha} U[k] \cup V_{\alpha} W[k] \right)
  	\setminus \large( V_{\alpha} (U \cap W)[k] 
  	\setminus V_{\alpha} (U \cap W)[2k+1] \large). 
  \]
\end{CorRef}

\begin{proof}
The subcomplexes in Corollary \refC{MNIH} now have a description based entirely on the enlargements of covering complexes and their intersections.
\end{proof}

\begin{RemRef}{GHGB}
In general, building $S[n]$ in $K$ for a subcomplex $S$ and a number $n \ge 1$ requires building $x[n]$ for all $x$ in $S$ and taking the union.  So the time complexity of this operation is proportional to the size of $S$.  Therefore, depending on the geometry of $K$ and the choices of $U$ and $W$ the sizes of the lists the computer needs to work though while building the enlargements can be greatly reduced by working with the intersection $U \cap W$ and one single enlarging procedure.  This is particularly true when the size of the intersection $U \cap W$ is significantly smaller than the sizes of $U$ and $W$.  In the case of simplicial trees which we address later in section \refS{EX}, the size of $U \cap W$ can be arranged to be an order of magnitude smaller than either $U$ or $W$ in the Euclidean case and that order can be made arbitrarily smaller in the case of a tree, depending on the valence of the tree. 
\end{RemRef}

In the rest of the section we illustrate in the specific context of the algorithm \texttt{ProcessLowerStars} \cite{RWS:11} that the general theorems can be improved on in situations with specific $\alpha$ and specific geometry.

Suppose $K$ is a cubulated 2D rectangle $I \times J$ and $I = [a,d]$ for $d-a \ge 4$.  Suppose $b$ and $c$ are integers nested in $a < b < c < d$.  Then we have a decomposition of $K$ as $U \cup W$ with $U = [a,c] \times J$ and $W = [b,d] \times J$.  

\begin{ThmRef}{NEJ}
Applying \texttt{ProcessLowerStars} as $\alpha$ to the two sets $U$ and $W$ as above, allows some savings in size of processed lists by using a more efficient formula
  \[
  	V_{\alpha} (K) = \left( V_{\alpha} (U) \cup V_{\alpha} (W) \right) \setminus
  	 V_{\alpha} (U\cap W) \setminus V_{\alpha} (U \cap W)[1]). \]
\end{ThmRef}

Notice that this theorem does not follow from any of the general theorems before because no enlargement of the patches is required for the containment $V_{\alpha} (K) \subset V_{\alpha} (U) \cup V_{\alpha} (W)$.

\begin{proof}
The proof is based on a case-by-case analysis of \texttt{ProcessLowerStars} applied to stars of vertices near $[b,c] \times J$.  As before, we are merging two lists $V_{\alpha} (U)$ and $V_{\alpha} (W)$.  We know that in each complex there are pairs that involve the vertices in $b \times J$ that may be included in non-critical pairs in error. In a complementary fashion, the same can be said about the critical cells in $C_{\alpha} (U)$ and $C_{\alpha} (W)$.  We know the errors can happen in this specific case of $\alpha$ because of a 2D counterexample similar to Example \refE{YIH}.  In much the same way the errors happen at this boundary because that is where incomplete information about their stars in $K$ is used within $U$ and $W$.  Let's call the pairs included in error ``mistakes''.  We try to identify the mistakes in the merged list $V_{\alpha} (U) \cup V_{\alpha} (W)$.  Are they guaranteed to be inside $V_{\alpha} (U \cap W)$?  This turns out to be true in this case but is not expected for
general uniformly local $\alpha$ or more general geometric decompositions.  Next, suppose a pair in $V_{\alpha} (U) \cup V_{\alpha} (W)$ is a mistake, then we know that it is not in $V_{\alpha} ((U \cap W)[1])$ because now both cells are in a star entirely contained in $V_{\alpha} ((U \cap W)[1])$.  This guarantees that the difference $ V_{\alpha} (U \cap W) \setminus V_{\alpha} ((U \cap W)[1])$ consists precisely of the mistakes in $V_{\alpha} (U) \cup V_{\alpha} (W)$ and completes the proof.
\end{proof}
 
\SecRef{Distributed \texttt{ProcessLowerStars} algorithm on a 2D digital image}{PA}

We apply Theorem \refT{NEJ} to parallelize an algorithm of type $\alpha$ on grayscale 2D images, in particular the algorithm \texttt{ProcessLowerStars} from \cite{RWS:11}.  

These digital images are modeled by cubical complexes with a grid-like structure, which makes it a little easier to work with than other cubical complexes, as we can index our patches with horizontal and vertical components, which will be outlined more clearly in a little bit.  Our formula for merging the vector field together only makes use of two patches.  In some cases, in particular, with a picture with many pixels and a computer program having a constraint on the number of pixels it can process with the algorithm, it is very likely that one would end up with more than two patches.  The grid-like structure that these cubical complexes possess will allow us to apply our formula to two pieces at a time.  First, the image will be split up into patches, which will partition the cubical complex. How many patches we end up with will be determined by the number of pixels in the image and the constraint we have on the number of pixels our computer program can process. After taking the star of each of the patches in both the vertical and horizontal directions separately, we will be left with new patches where there will be overlap, both vertically and horizontally, between adjacent patches.  This will allow us to apply our formula.  

After the partitioning into patches and taking the stars in both directions of all the patches, we apply the \texttt{ProcessLowerStars} algorithm to each patch and all of the overlaps of the adjacent patches, including the overlaps of the overlaps. The algorithm will give the discrete vector field in each of the patches and all of the overlaps of adjacent patches. We then start applying our formula, starting in upper left hand side of the our complex, with the first two starred patches.  This will give the discrete vector field of the union of these two starred patches.  We then proceed by applying our formula to this new bigger patch with the adjacent starred patch to its right.  We continue until we reach the end of this horizontal strip.  We then move down vertically to the next row of starred patches and go through the same procedure.  We do this for every horizontal strip of starred patches and obtain the discrete vector field for each horizontal strip.

We then move to the overlaps of the horizontal strips and use the same procedure as the horizontal strips themselves, starting from the left hand side, applying our formula as we move to the right end of the strip.  We will eventually end up with the discrete vector field of the overlaps of the horizontal strips.  We will need the discrete vector field of the overlaps of the horizontal strips in the next part of our algorithm when we start applying our formula vertically

We now have the discrete vector field of all of the horizontal strips and their overlaps, so we can now proceed applying our formula in the vertical direction, starting with the top most horizontal strip and applying our formula with the strip directly below it, making use of the discrete vector field of the overlaps that we computed in the previous steps.  We continue downward until we have the discrete vector field of one big patch and the horizontal strip below it, along with the discrete vector field of their overlap.  We apply our formula to these last two patches and end up with the discrete vector field of the whole image.

We make this precise with a pseudo-code for the distributed \texttt{ProcessLowerStars} algorithm. We assume that the constraint of our computer program in applying the \texttt{ProcessLowerStars} algorithm on a digital image is $N$ pixels. Assume that a digital image has $D$ pixels where $D>N$.  We break the image up into $n$ disjoint patches such that $n > {D}/{N}$, where patch $P_{i,j}^\#$ is a subcomplex and corresponds to the patch in the $ith$ column and $jth$ row for all $i=1,\cdots,m$ and $j=1, \cdots, \ell$ such that $m\cdot \ell=n$.  

Let $P_{i,j}^\#[a,b]$ be the subcomplex that is the $a$-neighborhood of $P_{i,j}^\#$ in the horizontal direction together with the $b$-neighborhood of $P_{i,j}^\#$ in the vertical direction.  We take our disjoint decomposition $P_{i,j}^\#$ and enlarge each patch by $1$-neighborhood in each direction and and call it $P_{i,j}$. This gives us $P_{i,j}=P_{i,j}^\#[1,1]$ where each patch now overlap and the overlaps contain at least one 2-cell. Let
\begin{enumerate}
\item $P_{i,j}^*=\displaystyle \bigcup_{k=1}^i P_{k,j}$ (the union of patches moving across a horizontal strip)
\item  $(P_{i,j} \cap P_{i, j+1})^*= \displaystyle \bigcup_{k=1}^i (P_{k,j} \cap P_{k, j+1})$ (the union of patches moving horizontally across the intersection of patches that are vertically adjacent to each other, i.e. moving horizontally along the intersection of strips)
\item $(P_{i,j}[0,1] \cap P_{i, j+1}[0,1])^*=\displaystyle \bigcup_{k=1}^i (P_{k,j}[0,1] \cap P_{k, j+1}[0,1])$ (the union of patches moving horizontally across the intersection of $1$-neighborhoods in the vertical direction of the patches that are vertically adjacent to each other) 
\end{enumerate}

\begin{figure}[h!]
      \centering
       \includegraphics[width=0.95\linewidth]{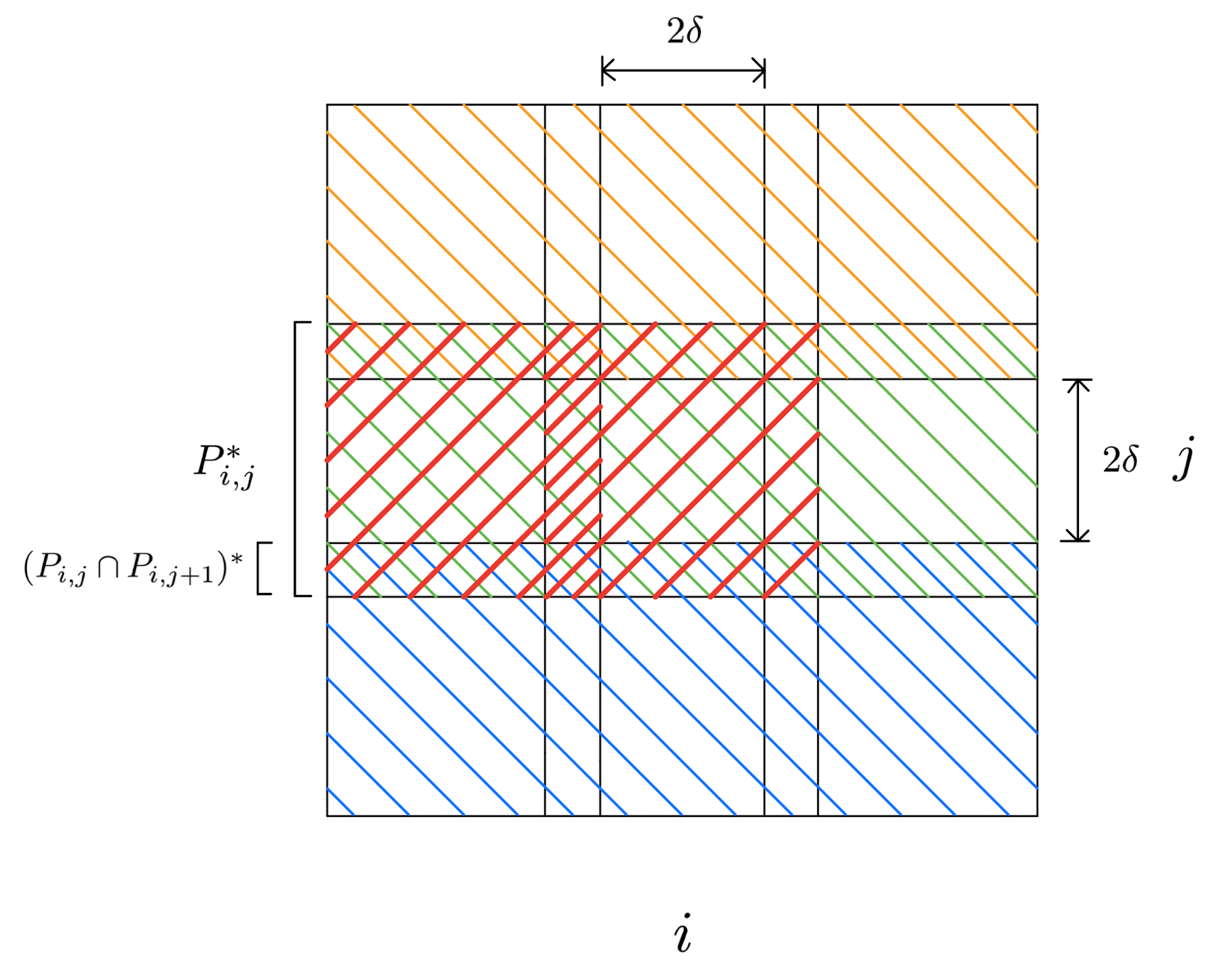} 
       \caption{An illustration of the planar grid-like decomposition.}
      \label{grid}
   \end{figure}
   
   Figure \ref{grid} is used to illustrate this notation.   The figure is intended to show a part of the decomposition of the 2D plane centered on the middle patch denoted $P_{i,j}$.  The orange, green, and blue overlapping rows are the unions of such patches across the values of the horizontal parameter $i$.  One useful observation is that there are only double overlaps of the rows so that there is always at least $2\delta$ clearance between non-adjacent rows such as the orange and the blue, for a fixed value of $\delta$.  This guarantees that $\delta$-enlargements of this covering of the plane by rows have the same nerve as the rows themselves. The nerve is isomorphic to a triangulation of the real line.  The figure also shows the set $P_{i,j}^*$, which is shaded red, and the union of intersections of patches $(P_{i,j} \cap P_{i, j+1})^*$.  It is clear that the finer covering of the plane by individual patches has the same property as above: its $\delta$-enlargement has the same nerve as the covering itself.

\begin{algorithm}
\caption{Distributed \texttt{ProcessLowerStars} on a 2D Digital Image}
\begin{algorithmic}
\STATE {\bfseries Inputs:} \begin{itemize}
\item $D$ digital image pixels
\item $g$ grayscale values on pixels
\end{itemize}
\STATE {\bfseries Outputs:} \begin{itemize}
\item $C$ critical cells
\item $V$ discrete vector field
\end{itemize}
\STATE {\bfseries Constraint:} \begin{itemize}
\item Memory$=N$ Pixels
\end{itemize}
\FOR {all $i$ and $j$} 
\STATE {Apply \texttt{ProcessLowerStars} algorithm to each of the following to obtain: 
\begin{itemize}
\item $P_{i,j} \rightarrow V_{i,j}$
\item $P_{i,j} \cap P_{i+1, j} \rightarrow V_{i,j}^{i+1}$
\item $P_{i,j} \cap P_{i, j+1} \rightarrow V_{i,j}^{j+1}$
\item $(P_{i,j} \cap P_{i, j+1}) \cap (P_{i+1,j} \cap P_{i+1,j+1}) \rightarrow V_{i,j}^{\cap}$
\item $P_{i,j}^*[1,0] \cap P_{i+1,j}[1,0] \rightarrow V_{i,j}^{i+1}[1,0]$
\item $(P_{i,j} \cap P_{i, j+1})^*[1,0] \cap (P_{i+1,j} \cap P_{i+1,j+1})[1,0] \rightarrow V_{i,j}^{\cap} [1,0]$
\item $P_{i,j}[0,1] \cap P_{i, j+1}[0,1] \rightarrow V_{i,j}^{j+1}[0,1]$
\item  $(P_{i,j}[0,1] \cap P_{i, j+1}[0,1]) \cap (P_{i+1,j}[0,1] \cap P_{i+1,j+1}[0,1]) \rightarrow V_{i,j}^{\cap}[0,1]$
\item $(P_{i,j}[0,1] \cap P_{i, j+1}[0,1])^*[1,0] \cap (P_{i+1,j}[0,1] \cap P_{i+1,j+1}[0,1])[1,0] \rightarrow V_{i,j}^{\cap}[1,1]$

\end{itemize}} \LONGCOMMENT{Preprocessing Step to obtain discrete vector field in all patches and their overlaps}
\ENDFOR
\FOR {$j=1, \cdots , \ell$}
\FOR {$i=1, \cdots, m-1$}
\STATE{\bfseries Update: $V_{i+1,j}^{'}=(V_{i,j}^{'} \cup V_{i+1,j}) \setminus (V_{i,j}^{i+1} \setminus V_{i,j}^{i+1}[1,0])$}\\
\LONGCOMMENT{Moving horizontally along strips}
\ENDFOR
\ENDFOR
\FOR {$j=1, \cdots , \ell$}
\FOR {$i=1, \cdots, m-1$}
\STATE{\bfseries Update: $(V_{i+1, j}^{j+1})^{'}=((V_{i, j}^{j+1})^{'} \cup V_{i+1, j}^{j+1})\setminus (V_{i,j}^{\cap} \setminus V_{i,j}^{\cap}[1,0])$}\\
\COMMENT{Moving horizontally along intersection of strips}
\ENDFOR
\ENDFOR
\FOR {$j=1, \cdots , \ell$}
\FOR {$i=1, \cdots, m-1$}
\STATE{\bfseries Update: $(V_{i+1,j}^{j+1}[0,1])^{'}=((V_{i,j}^{j+1}[0,1])^{'} \cup V_{i+1, j}^{j+1}[0,1]) \setminus (V_{i,j}^{\cap}[1,0] \setminus V_{i,j}^{\cap}[1,1])$}\\
\COMMENT{Moving horizontally along intersection of strips enlarged by 1 vertically}
\ENDFOR
\ENDFOR
\FOR {$j=1, \cdots , \ell-1$}
\STATE{\bfseries Update: $V_{j+1}^{'}=(V_{m,j}^{'} \cup V_{m,j+1}^{'}) \setminus ((V_{m,j}^{j+1})^{'} \setminus (V_{m,j}^{j+1}[0,1])^{'})$}\\
\COMMENT {Moving vertically down strips}
\ENDFOR
\end{algorithmic}
\end{algorithm}

\newpage

\SecRef{Generalization to a hierarchical tree-like decomposition}{EX}

Our main motivation for this work has been parallelization of the specific algorithm \texttt{ProcessLowerStars} on 2D images that we have done in sections \refS{FT} and \refS{PA}.  This algorithm has been proven to work in cubical cellular complexes based on standard cubical grids in 2D and 3D \cite{RWS:11}.  However, our theorems in section \refS{FT} have only general local constraints on the type of the algorithm and on the geometry of the regular cellular complex.  In this section we want to leverage some geometric material from \cite{BG} and \cite{GCBG}
to construct required antithetic coverings for grids in $n$D for all $n$ and, more generally, any subcomplex of a product of finite locally finite trees.  

It is known from a result of Dranishnikov \cite{Dranishnikov,Kasprowski}
that all metric spaces which satisfy a very weak and natural geometric condition called \textit{finite asymptotic dimension} (FAD) can be coarsely embedded in a finite product of locally finite trees with uniform distortion. This allows us to give a useful antithetic covering of any cellular complex which can be given an FAD metric with a universal bound on the size of all cells. 

As we observed before in Remark \refR{GHGB}, the importance of the case of a tree or a product of trees is also as an illustration of how crucial the improvements can be in passage from using the most general Theorem \refT{NEJCHCB} to using Corollary \refC{NEJCHCB}.

We saw in the previous section a worked-out example of use of antithetic decompositions in the case of a cubical grid in 2D.
Just as in that example, it is most natural to decompose a multi-parameter geometry according to projections to subsets in one chosen parameter. 
There should result an inductively defined decomposition of the entire cubical complex.
The general kind of parameter for our purposes is tree-based, with partial order, generalizing from the totally ordered real line.

A simplicial tree is a simplicial complex that is connected and has no cycles.  Recall also that a nerve of a collection of subsets of a given set is the simplicial complex with vertices which are the subsets and simplices corresponding to non-empty intersections of families of subsets.

\begin{DefRef}{HTLD}
  Suppose $\mathcal{U}$ is covering of a metric space.  We will denote by $\mathcal{U}[k]$ the covering by $k$-enlargements of members of $\mathcal{U}$ for a number $k \ge 0$.  
  We will say that a covering $\mathcal{U}$ is a \textit{tree-like decomposition with margin $k$} if $\mathcal{U}[k]$ is a simplicial tree.  Suppose each of the covering sets with the subspace metric also has a tree-like decomposition with margin $k$.  Then we say that the resulting covering by smaller sets is a \textit{hierarchical tree-like decomposition with margin $k$ and depth 2}.  Inductively, for a natural number $D$ one defines a \textit{hierarchical tree-like decomposition with margin $k$ and depth $D$}. We will refer to the sets that appear in such  hierarchical decomposition and are not unions of other sets as \textit{primary sets}. 
\end{DefRef}

The most useful hierarchical tree-like decomposition of depth $D$ can be obtained for any subset of the product of $D$ simplicial trees. The simplest case of this type of decomposition is when the trees are obtained as triangulations of the real line, which generalizes 2D decompositions as in Algorithm 2 to higher dimensions.

\begin{DefRef}{NMNYCF}
Let $T$ be a simplicial tree where each edge is given length 1 with the global metric induced as a path metric.  We fix a vertex $v_0$.  Given another vertex $v$ in $T$, we define the ``jet'' subset
$\J (v) = \{ t \in T \vert v \in [v_0, t) \}$.
Let $B(v,r)$ stand for the open metric ball of radius $r$ centered at $v$ and $S(v,r)$ stand for its boundary sphere.
We also define the subsets 
$\J (v,l) = \J (v) \cap B(v,l)$ 
for a positive number $l$, and the differences
$\J (v; l_1, l_2) = \J (v, l_2) \setminus B (v, l_1) \setminus S(v,l_1)$,
for $l_2 > l_1 > 0$.
\end{DefRef}

Given a number $r$ greater than 1, consider the collection of open subsets of $T$ consisting of the ball $B(v_0,2r)$ and the differences $\J (v; r-1,3r)$ where the vertices $v$ vary over $S(v_0, (2n-1)r)$ for arbitrary natural numbers $n$.  It is easy to see the following properties.
\begin{enumerate}
    \item This collection of subsets $\mathcal{V}$ is a covering of $T$.  Its nerve is itself in general a forest of trees where the vertices can be indexed by $v_0$ and the vertices $v \in S(v_0, (2n-1)r)$.  The edges are the pairs $(v,v')$ where $v' \in \J (v, 2(n+1)r)$.
    \item The diameter of each set in the covering is bounded by $6r$.
    \item The covering has a margin at least $r$, in the sense that the $r$-enlargement is again a tree isomorphic to the one in (1). 
\end{enumerate}
Clearly, intersecting any subset $X$ of $T$ with the produced covering gives a covering of $X$ with exactly the same three properties except possibly a forest of trees instead of a single tree.  This gives a hierarchical tree-like decomposition of $X$ with margin $r$ and depth $1$.

\begin{figure}[h!]
      \centering
       \includegraphics[width=0.70\linewidth]{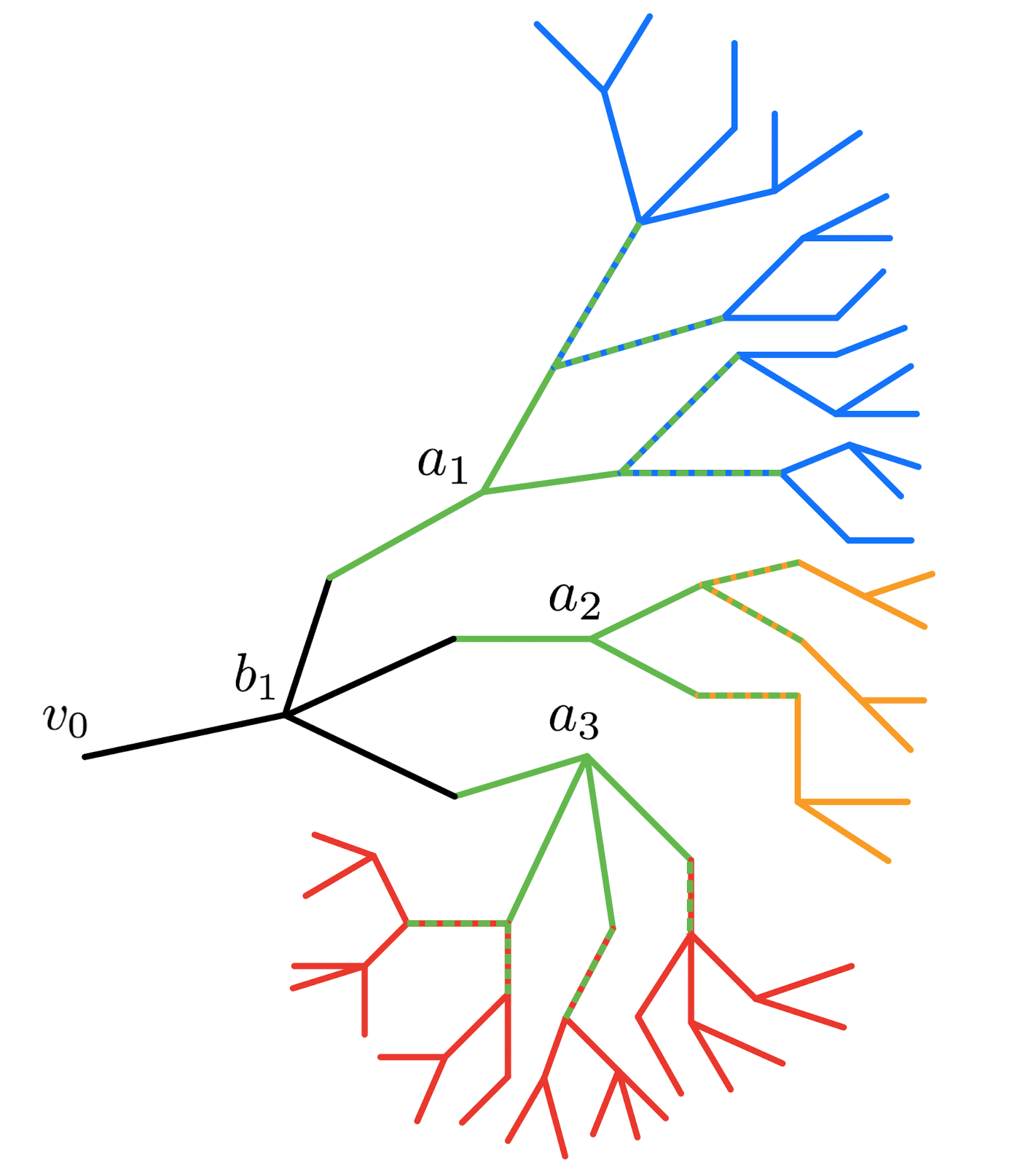} 
       \caption{An illustration of jet subsets $\J (v; l_1, l_2)$ and margins in a tree.}
      \label{tree}
   \end{figure}

Figure \ref{tree} is used to illustrate this notation.
In this figure, $v_0$ is the root of the tree.  Notice that the vertices labeled $a_1$, $a_2$, $a_3$ are all at distance 3 from $v_0$.  The blue, the orange, and the red subsets can be described as the jet subsets $J(a_i, 1, 4)$ for $1 \le i \le 3$.  They overlap with the green subset which itself can be described as the jet subset $J(b_1, 1, 4)$.  Notice that the bi-colored edges represent the overlaps between the different colored subsets.  This should illustrate the pattern of jet generated covering sets and their overlaps in the whole tree.  One more feature that can be seen in this figure is that with this particular set of choices of the vertices and the bounds $l_1$, $l_2$ there are only double overlaps.  There are no overlaps between jet subsets of higher multiplicities.  In other words, the nerve of such covering is 1-dimensional.  Now it can be easily seen in the picture that the cellular enlargement of all subsets by 1 unit still has a 1-dimensional nerve.  The corresponding enlargements are simply the jet subsets $J(a_i, 0, 5)$, for $1 \le i \le 3$, and $J(b_1, 0, 5)$.  This illustrates how this procedure with the specific choices we make produces a family of patches with margin 1. 

For a product $\Pi$ of $D$ trees, one starts by performing this construction in each of the factors, then takes the product of the covering families to create a covering of $\Pi$.  This is a hierarchical tree-like decomposition with margin $r$ and depth $D$.  Again, intersecting any subset $X$ of the product with the covering subsets is a hierarchical tree-like decomposition of $X$ with margin $r$ and depth $D$.

\begin{RemRef}{GASS}
Recall that a \textit{real tree} in geometry is a geodesic metric space where every triangle is a tripod.  Divergence behavior of geodesics in a real tree allows to generalize the above constructions in simplicial trees verbatim to real trees and their products.
\end{RemRef}

In the rest of this section, we show how to use a hierarchical tree-like decomposition similar to Algorithm 2.

Suppose our complex $K$ is a subcomplex of $\Pi$ with the product cubical structure.  In each factor we assume the simplicial tree structure. Distributing the computation, one simplifies the computation by applying an algorithm to slices which are reduced in size compared to the total complex. 

In our case this process will be performed inductively using the covering sets $\mathcal{V}$ in each tree coordinate.  
This will guarantee that all slices used in the computation form coverings with margin at least $r$.  We remind the reader that because of the assumptions both $\mathcal{V}$ and $\mathcal{V}[r]$ can be viewed themselves as vertices of forests of trees.

If $\pi_i$ is the projection onto the $i$-th coordinate, 
we have the slices $S_{i,V} = \pi_i^{-1} (V)$ for all $V$ in $\mathcal{V}_{i}$.  Since $\pi_i$ is a 1-Lipschitz function, the slices form a covering of $K$ with margin at least $r$.
We think of $\mathcal{V}_{i}$ as analogues of the intervals in Algorithm 2, so the following terminology is natural.

\begin{NotRef}{NBBBB}
Given a vector $x$ with $i$-th coordinate a set $V_i$ from $\mathcal{V}_i$, there are ``hypercubes''
\[
C_x = S_{1,V_1} \cap S_{2,V_2} \cap \ldots \cap S_{D, V_D}.
\]
\end{NotRef}

These hypercubes can serve as cells in a grid-like structure in $\Pi$.
Note that the empty set $\emptyset$ can be made a valid value for $V_{\ast}$.
Let $x$ be a vector as above with the property that if $\emptyset$ appears as a value then all subsequent coordinates must be $\emptyset$.  We will use $d(x)$ to denote the highest index for which the value is non-empty.  
Now the subsets of indices
$$X(s) = \{ x \in \prod \mathcal{V}_{i} \mid D-d(x)=s \}$$
and the corresponding coverings $\mathcal{C} (s)$ by $C_x$ for $x \in X(s)$
for increasing values of $s$ form the covering sets that generalize the rectangles and strips from Algorithm 2.
This is precisely the inductive structure that was leveraged in the 2D plane in the case of each tree $T_i$ being a simplicial line covered by the intervals.

\SecRef{Remarks on applications}{D}

There is now a large library of applications of discrete gradient fields through its use to simplify computations in topological data analysis but also direct applications to concrete problems.  Our treatment of its parallelization in this paper is very general.  To our knowledge all useful gradient fields that appear in the literature are uniformly local and so the methods of this paper can be applied to them.  

Just to give an example of gradient fields used for motion planning in robotics, we mention a couple of recent papers of the second author.  The problem in motion planning is to create an algorithm for a path a robot would take through the so-called Free Space of all allowable configurations, the configurations of the robot that avoid all obstructions in the physical environment.  The Free Space is really a parameter space of feasible robot positions.  One approach to modeling it is to sample the Free Space and build the corresponding simplicial approximation to it through some Vietoris-Rips complex built on the sample as vertices \cite{UEW}.  There are several ways to produce a discrete Morse function on a simplicial complex which restricts to given values on the vertices, or essentially equivalently a gradient field, some well very known such as \cite{King}.  All of them are uniformly local algorithms.  
One way this can be used for practical motion planning is as follows.  Convex polyhedral obstructions usually produce convex polyhedral exclusion zones complementary to the Free Space. One may chose a density estimator for the sample in Free Space, with values at every sample point, which is then possible to extend to a Morse function.  It is argued in \cite{IROS} that critical marginal points can act as ``lighthouses'' in planning a small finite collection of paths that provably contains an optimal path of smallest length which can then be easily extracted.

There are two kinds of applications of our formula in this setting.  The construction of the vector field can be distributed by considering an arbitrary or antithetic covering of the Free Space.  This allows one to construct the vector field that identifies the lighthouses in patches.  There is also the option of processing the data ``as needed'', for example by exploring an adjacent patch and the next lighthouse only when the robot approaches its border.

\SecRef{Discussion}{Disc}

There are plenty of comments in the literature that point out that by its nature \texttt{ProcessLowerStars} is embarrassingly parallelizable, cf. Robins et al. \cite{RWS:11} itself and, for another example, section 2 of Guylassy et al. \cite{Guylassy}. In the same way other uniformly local algorithms $\alpha$ are embarrassingly parallelizable.  This perspective of parallelization is not that useful unless $\alpha$ is very expensive to compute by itself.  A more urgent need is the type of distributed computation that builds partial vector fields on patches and combines them together.  Since discrete vector fields are analogues of smooth vector fields, it's not surprising that they are also very locally defined and are not expensive to compute in each locality.  These observations justify our distributed computation in this paper as the useful perspective on parallelization in this context.

Related to the point above, we wonder if there are useful discrete vector fields that are not uniformly local or even those that are uniformly local for $k \ge 2$ but not for $k=1$.  The authors are not aware of any in the literature.  Just to mention another example of a broadly used algorithm, it's easy to see that \texttt{ExtractRaw} of King et al. \cite{King} which is used to generate a discrete vector field from values of a Morse function on 0-cells is uniformly local for $k=1$.  If we are correct then our method in the paper applies to all known discrete vector fields.

Our future plans include extending Theorem \refT{NEJ} and Algorithm 2 to the 3D case where \texttt{ProcessLowerStars} has been proven to work \cite{RWS:11}.  We will also address similar custom efficiency improvements of general formulas for another algorithm \texttt{MorseReduce} due to Mischaikow and Nanda \cite{MiNa}.  This is a versatile algorithm using discrete Morse theory as preprocessing tool for persistent homology computations in topological data analysis.  The striking feature of \texttt{MorseReduce} is that it is dimension-independent and can be applied to very general regular cellular complexes.

Our results apply to infinite complexes and infinite coverings of complexes.  This comment might seem to have no practical implications, however we can use it to point out that our parallelization theorems allow a dynamic approach to processing the data. As in the robotics application described in the preceding section, one may not want to process all available patches at the same time but rather proceed one patch at a time depending on the need of a dynamic process such as planning a path.  In this case there is a need to incorporate the new patch into the pre-existing framework.  There might be infinitely many possible patches in the agnostic planning process.  The theorems can be used to extend the definition of the discrete vector field to each successive patch as many times as needed.

 Finally, we would like to contrast our theorems with comparable parallelization algorithms \cite{Guylassy08} and \cite{Shivashankar}.  We are grateful to the referee who pointed out these algorithms to us.  Both papers deal with discrete models approximated by smooth gradient vector fields and geometrically parallelize the computation of the associated Morse-Smale complexes.  The goals of these authors are very much akin to ours.  In fact, we refer the reader to additional motivation in Related Work sections in both papers. The main distinguishing feature of our theorems from section \refS{FT} is their generality.  They apply to any uniformly local algorithm on any regular cellular complex of any dimension, while the results of the referenced papers are more specific, geared toward computation of Morse-Smale complexes associated to gradient flows on lower dimensional manifolds. We don't say these are special cases of the general theorems.  They are certainly leaner and more efficient algorithms for that specific task, much like our Theorem \refT{NEJ} is not a special case of Theorem \refT{NEJCHCB} and its corollaries.


\begin{thebibliography}{99}

\bibitem{GCBG}
{G. Carlsson and B. Goldfarb},
\textit{The integral $K$-theoretic Novikov conjecture for groups with finite asymptotic dimension},
Invent. Math. \textbf{157} (2004), 405--418. 

\bibitem{DRS:15}
{O. Delgado-Friedrichs, V. Robins, and A. Sheppard}, {\em Skeletonization and partitioning of digital images using discrete Morse theory}, IEEE Trans. Pattern Anal. Mach. Intell. \textbf{37} (2015), 654--666.

\bibitem{Dranishnikov}
A. Dranishnikov,
\textit{On hypersphericity of manifolds with finite asymptotic dimension}, Trans.
Amer. Math. Soc. \textbf{355} (2003), 155--167.

\bibitem{F:02}
{R. Forman}, 
\textit{A user’s guide to discrete Morse theory}, {S{\'e}m. Lothar. Combin}, \textbf{48}, (2002), 35pp.

\bibitem{BG}
B. Goldfarb,
\textit{Singular persistent homology with geometrically parallelizable computation}, Topol. Proc. \textbf{55} (2020), 273--294.

\bibitem{Guylassy08}
A. Gyulassy, P.-T. Bremer, B. Hatmann and V. Pascucci, 
{\em A Practical Approach to Morse-Smale Complex Computation: Scalability and Generality}, in IEEE Trans. Vis. Comput. Graph. \textbf{6} (2008), 1619--1626.

\bibitem{Guylassy}
A. Gyulassy, P.-T. Bremer, and V. Pascucci,
{\em Computing Morse-Smale Complexes with Accurate Geometry}, IEEE Trans. Vis. Comput. Graph. \textbf{18} (2012), 2014--2022.

\bibitem{Kasprowski}
D. Kasprowski,
{\em Coarse embeddings into products of trees},
to appear in Kyoto Journal of Mathematics, \texttt{arXiv:1810.13361}.

\bibitem{King}
H. King, K. Knudson, and N. Mramor,
\textit{Generating discrete Morse functions from point data},
{Experimental Mathematics}
\textbf{14} (2005), 435--444.

\bibitem{Knudson}
{K.P. Knudson}, 
\textit{Morse Theory: Smooth and Discrete}, World Scientific, 2015.

\bibitem{MiNa}
{K. Mischaikow and V. Nanda},
{\em Morse Theory for Filtrations and Efficient Computation of Persistent Homology},
Discr. Comput. Geom. \textbf{50} (2013), 330--353.

\bibitem{RWS:11}
{V. Robins, P.J. Wood, and A. Sheppard}, {\em Theory and algorithms for constructing discrete Morse complexes from grayscale digital images}, IEEE Trans. Pattern Anal. Mach. Intell. \textbf{33} (2011), 1646--1658.

\bibitem{Shivashankar}
N. Shivashankar, Senthilnathan M, and V. Natarajan, 
{\em Parallel Computation of 2D Morse-Smale Complexes}, in IEEE Trans. Vis. Comput. Graph. \textbf{18} (2012), 1757--1770.

\bibitem{UEW}
{A. Upadhyay, W. Wang, and C. Ekenna},
\textit{Approximating cfree space topology by constructing Vietoris-Rips complex} in Proceedings of 2019 IEEE/RSJ International Conference on Intelligent Robots and Systems (IROS 2019), 2019, 2517--2523.

\bibitem{IROS}
{A. Upadhyay, B. Goldfarb, and C. Ekenna},
\textit{A topological approach to finding coarsely diverse paths}, in Proceedings of 2021 IEEE/RSJ International Conference on Intelligent Robots and Systems (IROS 2021).
 
\end{thebibliography}
\end{document}